\RequirePackage{fix-cm}
\documentclass[smallextended]{svjour3}       
\smartqed  
\usepackage{graphicx}

\usepackage{amsmath}
\usepackage{epsfig}
\usepackage{colordvi}
\usepackage{graphics}
\usepackage{color}

\usepackage{cite}
\usepackage{rotating}
\usepackage{booktabs}

 \usepackage{mathptmx}      
 \journalname{Journal of Elasticity}

\newfont{\tenbfit}{cmbx10}%

\newfont{\tenbbb}{msbm10}%
\newfont{\svnbbb}{msbm8}%

\newcommand{\lj}{\mbox{$[\kern-0.1478125em[$}}
\newcommand{\rj}{\mbox{$]\kern-0.1478125em]$}}
\newcommand{\la}{\mbox{$\langle\kern-0.2325em\langle$}}
\newcommand{\ra}{\mbox{$\rangle\kern-0.2325em\rangle$}}
\newcommand{\Blj}{\mbox{$\Big[\kern-0.275em\Big[$}}
\newcommand{\Brj}{\mbox{$\Big]\kern-0.275em\Big]$}}
\newcommand{\Bla}{\mbox{$\Big\langle\kern-0.425em\Big\langle$}}
\newcommand{\Bra}{\mbox{$\Big\rangle\kern-0.425em\Big\rangle$}}

\begin{document}

\title{Translation of Michael Sadowsky's paper ``An elementary proof for the existence of a developable \textsc{M\"obius} band and the attribution of the geometric problem to a variational problem"\footnote{Citations of this translation should refer also to Sadowsky's original paper, as cited in the Abstract.}
}

\titlerunning{An elementary proof for the existence of a developable \textsc{M\"obius} band}        

\author{Denis F.\ Hinz         \and
        Eliot Fried 
}


\institute{Denis F. Hinz \at
              Mathematical Soft Matter Unit\\
              Okinawa Institute of Science and Technology \\
              Okinawa, Japan 904-0495\\
              \email{dfhinz@gmail.com}        
           \and
           Eliot Fried \at
              Mathematical Soft Matter Unit\\
              Okinawa Institute of Science and Technology, \\
              Okinawa, Japan 904-0495\\
              \email{eliot.fried@oist.jp}
}

\date{Received: date / Accepted: date}

\maketitle

\begin{abstract}
This article is a translation of Michael Sadowsky's original paper ``Ein elementarer Beweis f\"ur die Existenz eines abwickelbaren \textsc{M\"obius}schen Bandes und die Zur\"uckf\"uhrung des geometrischen Problems auf ein Variationsproblem." which appeared in \emph{Sitzungsberichte der Preussischen Akademie der Wissenschaften}, \emph{physikal\-isch-mathematische Klasse}, 17.\ Juli 1930. --- Mitteilung vom 26.\ Juni, 412--415. Published on September 12, 1930. 
\end{abstract}


\section*{Translation of the original paper}

\textsc{M\"obius}~\cite{Moebius1865} illustrated the band bearing his name by describing how one may be constructed by bending a rectangular sheet of paper. Subsequently, it has been asked whether this construction can be achieved solely as a consequence of the compliance of the sheet in bending or whether stretching is also required. In other words --- the question has been raised as to whether the \textsc{M\"obius} band is developable in the strict sense.

\textsc{M\"obius} himself did not address this question, since it was irrelevant for his purposes. The developability of his band has been challenged by many, including, it is rumored, {\sc{H.\ A.\ Schwarz}}.
%
%
In the present work, the existence of a developable band is established on elementary geometric grounds. 

\begin{figure}[!t]
\includegraphics[width=0.75\textwidth]{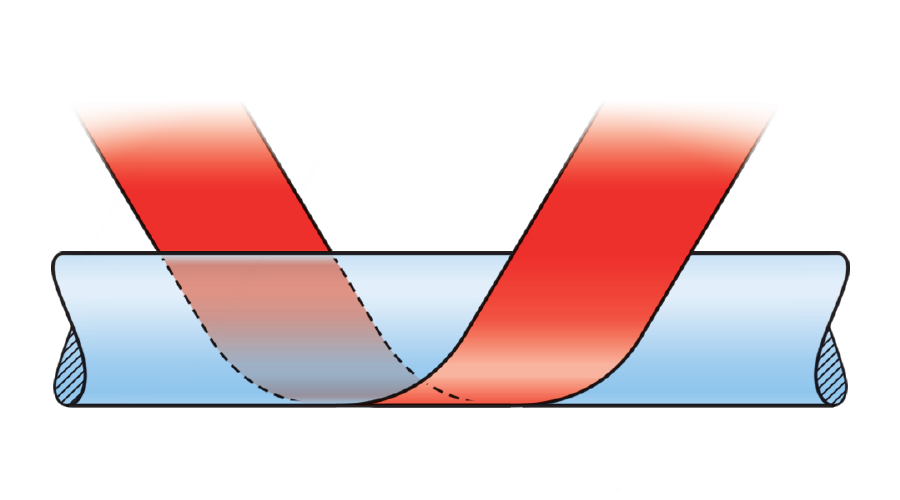}
 \caption{Adaptation of Figure~1 from the original version of the paper.}
 \label{fig01}
\end{figure}

Imagine an elongated rectangle of flexible but completely inextensible paper. 
Further, imagine a rigid cylindrical rod with circular cross-section. Construct two parallel planes tangent to the rod. A bendable but inextensible paper strip may then be positioned to lie within one of these planes, wrap halfway around the rod, and lie within the remaining plane. See Figure~\ref{fig01} for a depiction of the described arrangement.

\begin{figure}[!t]
\includegraphics[width=0.75\textwidth]{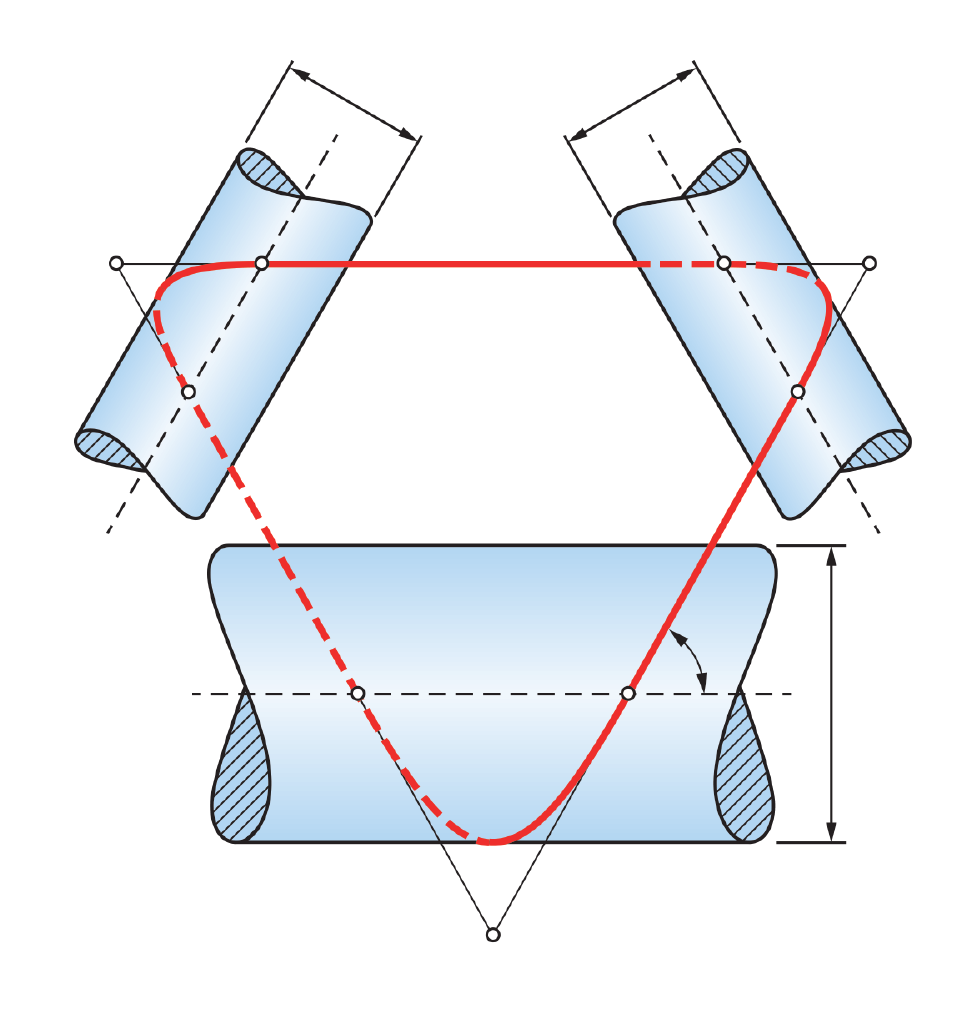}
\put(-121,22){\large$P$}
\put(-155,91){\large$A^\prime$}
\put(-98,91){\large$B^\prime$}
\put(-33,87){\large$2d$}
\put(-70,102){\large$60^\circ$}
\put(-58,165){\large$C^\prime$}
\put(-21,203){\large$R$}
\put(-74,191){\large$D^\prime$}
\put(-193,205){\large$E^\prime$}
\put(-235,203){\large$Q$}
\put(-198,165){\large$F^\prime$}
\put(-162,249){\begin{rotate}{-30}\large$d$\end{rotate}}
\put(-90.5,246){\begin{rotate}{30}\large$d$\end{rotate}}

 \caption{Adaptation of Figure~2 from the original version of the paper.}
 \label{fig02}
\end{figure}
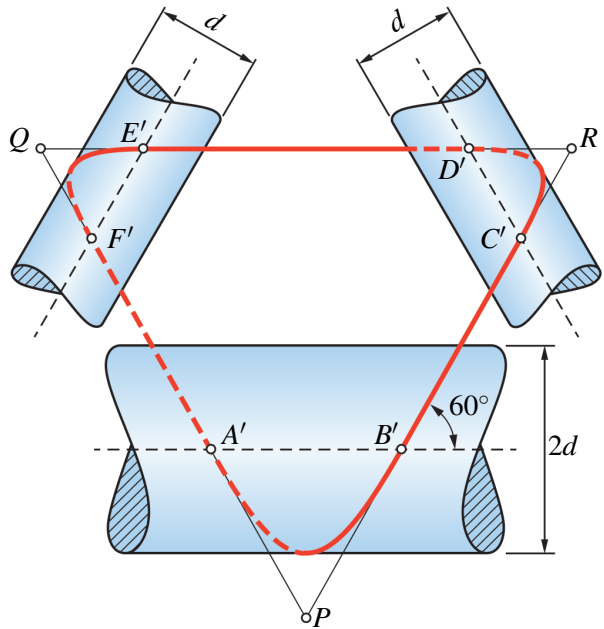

Consider now three cylindrical rods with circular cross-sections, one of which has diameter equal to the sum of the diameters of the remaining two. The rectangular strip of paper may then be threaded between the rods to form a \textsc{M\"obius} band. The resulting surface consists of three planar sections and three semicylindrical sections and, thus, is developable. An illustration of the construction for rods of diameters $d$, $d$, and $2d$ is provided in Figure~\ref{fig02}.
The elementary connections between the particular angles and lengths needed to form such a band are not of interest here, since only the existence of the band matters. 

The foregoing descriptive construction of a developable \textsc{M\"obius} band admits an equivalent analytical construction. If attention is restricted to the previously considered simple case in which the diameters of the rods are $d$, $d$, and $2d$, the projections onto the drawing plane of the midlines of the rectangular portions of the band lie on the edges of an equilateral triangle with side length $L$, and the axes of all three rods are parallel to the drawing plane. Let $A$, $B$, $C$, $D$, $E$, and $F$ be points on the midline of the \textsc{M\"obius} band and let $A^\prime$, $B^\prime$, $C^\prime$, $D^\prime$, $E^\prime$, and $F^\prime$ be their projections onto the drawing plane. Let $AB$, $BC$, etc.\ denote the distances measured along the mid-line from point $A$ to $B$, $B$ to $C$, etc.\ and let $A^\prime B^\prime$, $B^\prime C^\prime$, etc.\ denote the lengths of the rectilinear connections between $A^\prime$ and $B^\prime$, $B^\prime$ to $C^\prime$, etc.\ in the drawing plane. Whereas the sections $AF$, $BC$, and $DE$ of the midline remain rectilinear, the sections $AB$, $CD$, and $EF$ are twisted along helical curves. These helical curves intersect the generators of the cylinder with an angle of $60^{\circ}$, from which it follows that (Figure~\ref{fig02})
\begin{equation}\label{eq:01}
\begin{array}{c}
\displaystyle
A^\prime B^\prime = AB \cos(60^{\circ}), \hspace{0.5cm} AB \cos(30^{\circ}) = \frac{1}{2} \pi 2 d = \pi d, \hspace{0.5cm} {\rm{etc.,}} 
\cr\noalign{\vskip8pt}
\displaystyle
AB = \frac{2\pi d}{\sqrt{3}}, \hspace{0.5cm} A^\prime B^\prime = \frac{\pi d}{\sqrt{3}},
\cr\noalign{\vskip8pt}
\displaystyle
FE = CD = \frac{\pi d}{\sqrt{3}}, \hspace{0.5cm} F^\prime E^\prime = C^\prime D^\prime = \frac{\pi d}{2\sqrt{3}},
\cr\noalign{\vskip8pt}
\displaystyle
PQ = QR = PR = L,
\cr\noalign{\vskip8pt}
\displaystyle
AF = BC = L - \frac{\sqrt{3} \pi d}{2},
\cr\noalign{\vskip8pt}
\displaystyle
ED = L- \frac{\pi d}{\sqrt{3}}.
\end{array}
\end{equation}
Moreover, the total length $l$ of the midline obeys $l = 3L$.

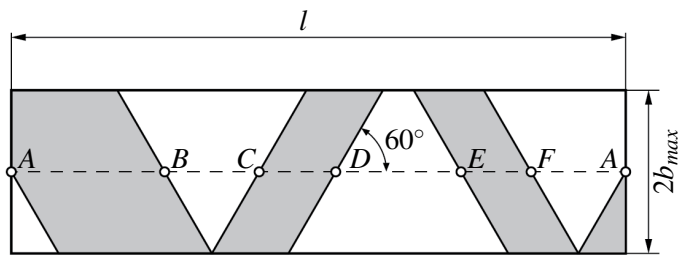
\begin{figure}[!t]
\includegraphics[width=0.75\textwidth]{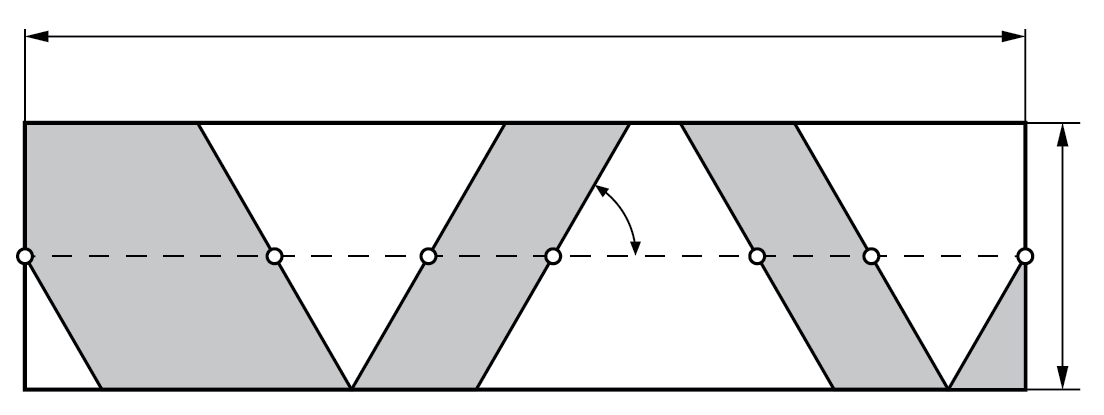}
\put(-245,38){\large$A$}
\put(-188,38){\large$B$}
\put(-163,38){\large$C$}
\put(-121,38){\large$D$}
\put(-108,45){\large$60^\circ$}
\put(-77,38){\large$E$}
\put(-51,38){\large$F$}
\put(-27,38){\large$A$}
\put(-140,90){\large$l$}
\put(0,28){\begin{rotate}{90}\large$2b_{max}$\end{rotate}}

 \caption{Adaptation of Figure~3 from the original version of the paper.}
 \label{fig03}
\end{figure}

Hereafter, $l$ and $d$ may be chosen arbitrarily (the latter within certain limits, cf.~\eqref{eq:02}) and, based on the computations above, the points $A$ to $F$ may be chosen to lie on the mid-line. These points also must be in the order $A$-$B$-$C$-$D$-$E$-$F$, which is the case for $3\sqrt{3} \pi d \leq 2 l$ (cf.~\eqref{eq:01}). Next, the generators of the cylindrically twisted parts can be marked on the band with an angle of $60^\circ$ (Figure~\ref{fig03}). The maximal band width $2b$ corresponding to one pair $(l,d)$  of values in the construction is determined through the generators, because the portions of the band, which are wound around \emph{different} cylinders may not overlap --- at most they may come into contact at the boundary of the band. Since $BC =  AF$ are the shortest of the unbent portions of the band, it follows that (Figure~\ref{fig03}) 

\begin{equation*}
\displaystyle
b \leq \frac{BC}{2}\tan(60^o), \quad {\rm{i.e.}} \quad b \leq \frac{l}{2\sqrt{3}}-\frac{3\pi d}{4}.
\end{equation*}

These inequalities identify the constraint
\begin{equation}\label{eq:02}
\displaystyle
d \leq \frac{2 l}{3 \sqrt{3} \pi}.
\end{equation}
For an infinitely narrow ($b=0$) \textsc{M\"obius} band, we may set
\begin{equation*}
\displaystyle
BC = AF = 0.
\end{equation*}
For such a band we also have
\begin{equation*}
\displaystyle
AB = \frac{4}{9}l, \quad FE=CD=\frac{2}{9}l, \quad DE = \frac{1}{9} l,
\end{equation*}
and (cf.\ the following considerations)
\begin{equation}\label{eq:03}
\displaystyle
\int \limits_0^l H^2 ds = \frac{15\pi^2}{l}
\end{equation}
as an upper bound for the integral of the general variational problem~\eqref{eq:04} for a band that is infinitesimally narrow with respect to its length.

The existence of a developable \textsc{M\"obius} band is therefore established by construction. 
Consider, in more detail, the shape of the band so obtained: 
The band is assembled from pieces of individual flat and cylindrical surfaces; the actual band constructed with a piece of paper has a different smoothly curved shape. We may arrive at this form if we take into account that the band constructed out of planes and cylinder surface pieces may remain in equilibrium only under the influence of boundary and surface forces, respectively, along with moments. This follows, for example, on recognizing that the curvature of the band as constructed is discontinuous, and, moreover, that moments, as is well known, are associated with curvature. However, if freed from the influence of these forces the constructed band will deform and find an equilibrium configuration corresponding to a minimum of the internal energy associated with elastic deformation. The elastic bending energy density is proportional to the surface integral of the squared mean curvature, granted that the surface is developable. For the situation under consideration, the following observation holds: As initial shape we had a strictly developable band. Corresponding to the provision that the band is inextensible, the deformation is a pure deflection --- without stretching. We thus conclude that any configuration resulting from such deformation must be strictly developable. In general, the shape of a the \textsc{M\"obius} band corresponds to a minimum of the bending energy. The determination of this ideal form is thus attributed to the variational problem
\begin{equation}\label{eq:04}
\displaystyle
\iint  H^2 dF = {\rm{Min}},
\end{equation}
the integration being over the complete surface, where $H$ is the mean curvature of the surface and the following conditions apply: the strip forming the band surface has rectangular shape and the band is both developable and possesses the proper one-sided spatial connectivity.

For a infinitely narrow band, that is a band with width $b$ infinitesimally small compared to its length $l$, the variational problem~\eqref{eq:04} reduces to
\begin{equation}\label{eq:05}
\displaystyle
\int \limits_0^l  \frac{(K^2+W^2)^2}{K^2} ds = {\rm{Min}},
\end{equation}
where $K$ and $W$ are the curvature and the torsion, respectively, of the mid-line of the band, and $s$ is its arc length. For a further \emph{exact} treatment of the problem, it seems that~\eqref{eq:05} is practically inappropriate, since the calculation becomes exceedingly complicated.
\\

\emph{Addendum during the review:} For further determination of the configuration of the band cf.\ the work of the author in the reports of the $3^{\rm{rd}}$ International Congress for Technical Mechanics (3.\ internationaler Kongress f.\ techn.\ Mechanik), Stockholm 1930.

\end{document}